\documentclass[a4paper,11pt,leqno,twoside]{article}
\usepackage{amsmath,amsthm,amssymb,yoshinori}
\numberwithin{equation}{section}

\layout

\newcommand{\limup}{\lim_{q \uparrow 1}}

\newcommand{\sumn}{\sum_{n \in \bZ\backslash\{0\}}}
\newcommand{\nn}{^{(\n)}}
\newcommand{\tp}[1]{^{(#1)}}
\newcommand{\B}{\mathrm{B}}
\newcommand{\tz}{\tilde{\zeta}}
\newcommand{\tG}{\widetilde{\Gamma}}
\newcommand{\tB}{\widetilde{B}}

\newcommand{\bsym}[1]{{\boldsymbol{#1}}}

\newcommand{\ua}{\uparrow}

\renewenvironment{proof}{%
\smallbreak
\noindent
\textbf{Proof\,.}}

\newenvironment{Proof}[1]{%
\smallbreak
\noindent
\textbf{Proof of {#1}\,.}}

\boldtitle{$q$-ANALOGUES OF THE BARNES MULTIPLE\\ ZETA FUNCTIONS}
\author{YOSHINORI YAMASAKI\\
\smallbreak\\
{\it Graduate School of Mathematics, Kyushu University,}\\
{\it Hakozaki Fukuoka 812-8581 Japan.}\\
{\it ma203032@math.kyushu-u.ac.jp}
}

\begin{document}

\setlength{\baselineskip}{16pt}
\maketitle

\begin{abstract}
 In this paper, we introduce $q$-analogues of the Barnes multiple zeta
 functions. We show that these functions can be extended meromorphically
 to the whole plane, and moreover, tend to the Barnes multiple zeta
 functions when $q\ua 1$ for all complex numbers.   
\begin{MSC}
 11M41, 11M35
\end{MSC}
\end{abstract}

\section{Introduction}

 The aim of the present paper is to introduce $q$-analogues of the
 Barnes multiple zeta function (\cite{Barnes1904});  
\[
 \z_r(s,z;\bsym{\om})
:=\sum_{n_1,\ldots,n_r\ge 0}(n_1\om_1+\cdots+n_r\om_r+z)^{-s} \qquad
(\Re(s)>r),  
\]
 where $\om_1,\ldots,\om_r$ are complex parameters which lie on some
 half plane. We study an analytic continuation of the $q$-analogue of
 $\z_r(s,z;\bsym{\om})$. We determine especially, {\it true}
 $q$-analogues of the Barnes multiple zeta function when $\om_i=1$
 ($1\le i \le r$). Here, by a true $q$-analogue, we mean when the
 classical limit $q\ua 1$ of the $q$-analogue reproduces the original
 zeta function for {\it all} $s\in\bC$. Recall the Hurwitz zeta 
 function's case, that is, the case $r=1$. Let $0<q<1$ and
 $[z]_q:=(1-q^z)/(1-q)$ for $z\in\bC$. In \cite{KawagoeWakayamaYamasaki}
 (see also \cite{KanekoKurokawaWakayama2003}) we studied $q$-analogues
 of the Hurwitz zeta function $\z(s,z):=\sum^{\infty}_{n=0}(n+z)^{-s}$
 defined via the $q$-series with two complex variables $s$, $t\in\bC$;  
\[
 \tz_q(s,t,z):=\sum^{\infty}_{n=0}\frac{q^{(n+z)t}}{{[n+z]_q}^s} \qquad
 (\Re(t)>0). 
\]
 The function $\tz_q(s,t,z)$ is continued meromorphically to the whole
 $s$, $t$-plane. We obtained the necessary and sufficient condition for
 the variable $t\in\bC$ so that $\tz_q(s,t,z)$ is a true 
 $q$-analogue of $\z(s,z)$. Namely, these functions
 $\tz\nn_q(s,z):=\tz_q(s,s-\n,z)$ $(\n\in\bN)$ give true $q$-analogues
 of the Hurwitz zeta function among the functions of the form
 $\tz_q(s,\vp(s),z)$ where $\vp(s)$ is a meromorphic function on
 $\bC$. The main purpose is to generalize the results in
 \cite{KawagoeWakayamaYamasaki} to $r\ge 1$.      

 The plan of this paper is as follows.
 In \Sec{def}, we define the $q$-analogue $\z_{q,r}(s,t,z)$ of the
 Barnes multiple zeta function for $\om_i=1$ ($1\le i \le r$) and give
 the main theorem (\Thm{main}). 
 In \Sec{proof}, we first study an analytic continuation of the
 $q$-analogue $\z_{q,r}(s,t,z)$ and then prove the main theorem.   
 In \Sec{general}, we study the $q$-analogue
 $\z_{q,r}(s,t,z;\bsym{\om})$ of the multiple zeta functions for general
 parameters $\bsym{\om}$. Using the binomial theorem, we give an
 analytic continuation of $q$-analogues (\Prop{ac-Gzqr}). 
 In the appendix, we introduce a $q$-analogue $\tG_q(z)$ of the gamma
 function $\G(z)$ associated to the $q$-analogue $\tz_q(s,t,z)$ of the
 Hurwitz zeta function. We first observe fundamental properties of
 $\tG_q(z)$. The rest of the appendix is devoted to study $q$-analogues
 of the limit formula of Lerch (\Prop{q-Lerch}) and the Gauss-Legendre 
 formula (\Prop{Gauss-Legendre}).       
 
 Throughout the paper, we assume $0<q<1$. We put 
 ${[n]_q}!:=[n]_q[n-1]_q\cdots [1]_q$ for $n\in\bN$. Further, for 
 non-negative integers $m$ and $n$, we define the $q$-binomial
 coefficient ${m \brack n}_q$ by       
\[
 {m \brack n}_q:=\frac{(q;q)_{m}}{(q;q)_n(q;q)_{m-n}},
\]
 where $(a;q)_m:=\prod^{m-1}_{l=0}(1-aq^l)$ for $m\ge 1$ and
 $(a;q)_0:=1$. We denote the field of complex numbers, the ring of
 rational integers and the set of positive integers by $\bC$, $\bZ$ and
 $\bN$ respectively. Also, if $Q$ is a set, $Q_P$ stands for the set
 of all elements in $Q$ which satisfy the condition $P$.

\section{Definition of $q$-analogues and the main theorem}
\label{sec:def}

 Let $s,t\in\bC$ and $z\notin -\bZ_{\le 0}$. We study a $q$-analogue of
 the Barnes multiple zeta function   
\begin{align*}
 \z_r(s,z)
:&=\sum_{n_1,\ldots,n_r\ge 0}(n_1+\cdots+n_r+z)^{-s}
\intertext{defined by the following $q$-series;}
 \z_{q,r}(s,t,z)
:&=\sum_{n_1,\ldots,n_r\ge
0}\frac{q^{n_1t+n_2(t-1)+\cdots+n_r(t-r+1)}}{{[n_1+\cdots+n_r+z]_q}^s}.
\end{align*}
 The series $\z_{q,r}(s,t,z)$ converges absolutely for
 $\Re(t)>r-1$. When $r=1$, we put $\z_q(s,t,z):=\z_{q,1}(s,t,z)$. In
 view of the results in \cite{KawagoeWakayamaYamasaki}, 
 we put $\z\nn_{q,r}(s,z):=\z_{q,r}(s,s-\n,z)$ and 
 $\z\nn_q(s,z):=\z_{q}(s,s-\n,z)$ for $\n\in\bN$. The following theorem
 is the main result of this paper.   

\begin{thm}
\label{thm:main}
 Let $t=\vp(s)$ be a meromorphic function on $\bC$. Then the formula  
\[
 \limup \z_{q,r}(s,\vp(s),z)=\z_r(s,z) \qquad (s\in\bC)
\]
 holds if and only if the function $\vp(s)$ can be written as
 $\vp(s)=s-\n$ for some $\n\in\bN$.
\end{thm}

\begin{remark}
\label{rem:many}
 $(i)$\ By \Thm{main}, it is clear that the functions of the type
 $\sum_{\n\,:\,\textrm{finite}}a\nn_q(s,z)\z\nn_{q,r}(s,z)$ for some
 holomorphic functions $a\nn_q(s,z)$ satisfying
 $\limup \sum_{\n\,:\,\textrm{finite}}a\nn_q(s,z)=1$ are also true
 $q$-analogues of $\z_r(s,z)$. Note that the $q$-analogue of the
 Hurwitz zeta function discussed in \cite{KawagoeWakayamaYamasaki} is 
 given by $\tz\nn_q(s,z)=\z\nn_q(s,z)\times q^{z(s-\n)}$. 

 $(ii)$\ The $q$-analogue of the Hurwitz zeta function studied in
 \cite{Tsumura2001} is different from ours. It is not form of the
 ($q$-) Dirichlet series and, in fact, is needed an extra term
 (precisely, see \cite[Corollary\,2.4.]{KawagoeWakayamaYamasaki}). 
\end{remark}

 It is easy to see that $\z_r(s,z)$ is expressed as
\begin{equation}
\label{for:binom-zr}
 \z_r(s,z)
=\sum^{\infty}_{n=0}\binom{n+r-1}{r-1}(n+z)^{-s}.
\end{equation} 
 To obtain a similar expression for $\z_{q,r}(s,t,z)$, we need the 
 following lemma.  

\begin{lem}
\label{lem:q-binom}
$(i)$\ For $l, m\in\bZ_{\ge 0}$, it holds that
\begin{equation}
\label{for:q-binom1}
 \sum^{l}_{d=0}{m-1+d \brack m-1}_q{q^d}={m+l \brack m}_q.
\end{equation}

$(ii)$\ For $r\in\bN$, it holds that
\begin{equation}
\label{for:q-binom2}
\sum_{n_1,\ldots,n_r\ge 0 \atop
 n_1+\cdots+n_r=n}q^{n_1+2n_2+\cdots+rn_r}
=q^n{n+r-1 \brack r-1}_q. 
\end{equation}
\end{lem}
\begin{proof}
 The formula \eqref{for:q-binom1} is well-known (see,
 \cite{AndrewsAskeyRoy1999}, also
 \cite{KimotoKurokawaMatsumotoWakayama2005}). We show the formula
 \eqref{for:q-binom2} by induction on $r$. It is clear that 
 \eqref{for:q-binom2} holds for $r=1$. Suppose it holds for $r-1$. Then
 the left hand side of \eqref{for:q-binom2} is equal to    
\[
 \sum^n_{n_1=0}q^{n_1+\cdots +n_r}\sum_{n_2,\ldots,n_r\ge 0 \atop
 n_2+\cdots+n_r=n-n_1}q^{n_2+2n_3+\cdots+(r-1)n_r}
=q^n\sum^n_{n_1=0}{n_1+r-2 \brack r-2}_q q^{n_1}.
\]
 Using the formula \eqref{for:q-binom1} for $l=n$, $m=r-1$ and $d=n_1$,
 we obtain the desired formula.\qed
\end{proof}

\begin{prop}
\label{prop:q-binom-zqr}
 It holds that
\begin{equation}
\label{for:q-binom-zqr}
 \z_{q,r}(s,t,z)
=\sum^{\infty}_{n=0}{n+r-1 \brack r-1}_q\frac{q^{n(t-r+1)}}{{[n+z]_q}^s}. 
\end{equation}
\end{prop}
\begin{proof}
 It is easy to see that 
\begin{align*}
 \z_{q,r}(s,t,z)
&=\sum^{\infty}_{n=0}\sum_{n_1,\ldots,n_r\ge 0 \atop
 n_1+\cdots+n_r=n}\frac{q^{(t+1)(n_1+\cdots+n_r)-(n_1+2n_2+\cdots+rn_r)}}{{[n_1+\cdots+n_r+z]_q}^s}\\
&=\sum^{\infty}_{n=0}\frac{q^{(t+1)n}}{{[n+z]_q}^s}\sum_{n_1,\ldots,n_r\ge 0 \atop  n_1+\cdots+n_r=n}q^{-(n_1+2n_2+\cdots+rn_r)}.
\end{align*}
 Substituting $q^{-1}$ for $q$ into \eqref{for:q-binom2} yields  
\[
 \sum_{n_1,\ldots,n_r\ge 0 \atop
 n_1+\cdots+n_r=n}q^{-(n_1+2n_2+\cdots+rn_r)}
=q^{-n}{n+r-1 \brack r-1}_{q^{-1}}=q^{-nr}{n+r-1 \brack r-1}_q.
\]
 Hence we obtain the formula \eqref{for:q-binom-zqr}.\qed
\end{proof}

\section{Proof of the main theorem}
\label{sec:proof}

 In this section, we give a proof of \Thm{main}. We first provide
 analytic continuations of $\z_r(s,z)$ with respect to $s$ (see
 \cite{SeoChoiGangOk1993}) and study of $\z_{q,r}(s,t,z)$ with respect
 to $t$. Since we have the following ladder relations 
\begin{align}
 \z_r(s,z)
&=\z_r(s,z+1)+\z_{r-1}(s,z),
\nonumber\\
 \z_{q,r}(s,t,z)
&= q^{t-r+1}\z_{q,r}(s,t,z+1)+\z_{q,r-1}(s,t,z),
\label{for:ladder zqr}
\end{align}   
 it is sufficient to study the analytic continuation when
 $\Re(z)>0$. Here we understand $\z_0(s,z)=z^{-s}$ and
 $\z_{q,0}(s,t,z)={[z]_q}^{-s}$.
 
\subsection{An analytic continuation of $\z_r(s,z)$}

 For each $l\in\bZ_{\ge 0}$, we put
 $(x)_l:=x(x+1)\cdots (x+l-1)=\G(x+l)/\G(x)$. Then $(x)_l$ can be
 written as $(x)_l=\sum^{l}_{j=0}s(l,j)x^j$ where $s(l,j)$ is the
 Stirling number of the first kind. Hence we have       
\[
 \binom{n+r-1}{r-1}
=\frac{(n)_{r}}{n(r-1)!}=\frac{1}{(r-1)!}\sum^{r}_{j=0}s(l,j)n^{j-1}
=\sum^{r-1}_{l=0}P^l_r(z)(n+z)^l,
\]
 where $P^l_r(z)$ ($0\le l \le r-1$) is a polynomial in $z$ defined by  
\[
 P^l_r(z):=\frac{1}{(r-1)!}\sum^{r-1}_{j=l}\binom{j}{l}s(r,j+1)(-z)^{j-l}. 
\]
 Thus, we have by \eqref{for:binom-zr} 
\begin{equation}
\label{for:zr-z}
 \z_r(s,z)=\sum^{r-1}_{l=0}P^l_r(z)\z(s-l,z).
\end{equation}
 Recall also the Euler-Maclaurin summation formula (see, e.g.,
 \cite[p.\,619]{AndrewsAskeyRoy1999})\,:\,For $a, b\in\bZ$ satisfying
 $a<b$, a $C^{\infty}$-function $f(x)$ on $[a,\infty)$, and an arbitrary integer
 $M \ge 0$, we have          
\begin{multline}
\label{for:EM} 
 \qquad \sum^b_{n=a}f(n)
=\int^b_a f(x)dx+\frac{1}{2}(f(a)+f(b))\\
+\sum^M_{k=1}\frac{B_{k+1}}{(k+1)!}(f^{(k)}(b)-f^{(k)}(a))   
-\frac{(-1)^{M+1}}{(M+1)!}\int^b_a\tB_{M+1}(x)f^{(M+1)}(x)dx,\qquad
\end{multline}
 where $B_k$ is the Bernoulli number and $\tB_{k}(x)$ is the periodic
 Bernoulli polynomial defined by $\tB_{k}(x)=B_{k}(x-\Gauss{x})$ with
 $\Gauss{x}$ being the largest integer not exceeding $x$. Putting
 $f(x):=(x+z)^{-s}$, we obtain 
\begin{multline}
\label{for:EM-z}
 \qquad \qquad \z(s,z)
=\frac{1}{s-1}z^{-s+1}+\frac{1}{2}z^{-s}+\sum^M_{k=1}\frac{B_{k+1}}{(k+1)!}(s)_kz^{-s-k}\\
-\frac{(s)_{M+1}}{(M+1)!}\int^{\infty}_0\tB_{M+1}(x)(x+z)^{-s-M-1}dx.\qquad\qquad\qquad
\end{multline}
 Since $\Re(z)>0$, the equation \eqref{for:EM-z} gives an analytic
 continuation of the Hurwitz zeta function $\z(s,z)$ to the region
 $\Re(s)>-M$. Therefore, by \eqref{for:zr-z} and \eqref{for:EM-z}, we
 obtain the following    

\begin{prop}
\label{prop:zr}
 For any integers $M_l\ge 0$ $(0\le l \le r-1)$, we have
\begin{multline*}
 \z_r(s,z)
=\sum^{r-1}_{l=0}\frac{P^l_r(z)}{s-l-1}z^{-s+l+1}+\frac{1}{2}\sum^{r-1}_{l=0}P^l_r(z)z^{-s+l}
+\sum^{r-1}_{l=0}P^l_r(z)\sum^{M_l}_{k_l=1}\frac{B_{k_l+1}}{(k_l+1)!}(s-l)_{k_l}z^{-s+l-k_l}\\
-\sum^{r-1}_{l=0}\frac{P^l_r(z)(s-l)_{M_l+1}}{(M_l+1)!}\int^{\infty}_0\tB_{M_l+1}(x)(x+z)^{-s+l-M_l-1}dx.
\end{multline*}
 This gives an analytic continuation of $\z_r(s,z)$ to the region
 $\Re(s)>M$ where $M:=\max\{-M_l+l\,|\,0\le l\le r-1\}$.\qed 
\end{prop}

\subsection{An analytic continuation of $\z_{q,r}(s,t,z)$}

 It is easy to see that
\begin{align*}
 {n+r-1 \brack r-1}_q
&=\frac{1}{{[r-1]_q}!}\prod^{r-1}_{j=1}\frac{1-q^{n+z}+q^{n+z}-q^{n+j}}{1-q}\\
&=\frac{1}{{[r-1]_q}!}\prod^{r-1}_{j=1}\Bigl([n+z]_q-q^{n+j}[z-j]_q\Bigr)
 =\sum^{r-1}_{l=0}q^{n(r-1-l)}P^l_{q,r}(z){[n+z]_q}^l,
\end{align*}
 where $P^l_{q,r}(z)$ ($0\le l \le r-1$) is a function of $z$ defined
 by      
\[
 P^{l}_{q,r}(z):= \frac{(-1)^{r-1-l}}{{[r-1]_q}!}\sum_{1\le
 m_1<\cdots<m_{r-1-l\le r-1}}q^{m_1+\cdots+m_{r-1-l}}[z-m_1]_q \cdots [z-m_{r-1-l}]_q 
\]
 for $0\le l\le r-2$ and $P^{r-1}_{q,r}(z):=1/{{[r-1]_q}!}$. Therefore
 we have by \eqref{for:q-binom-zqr} 
\begin{equation}
\label{for:zqr-zq}
 \z_{q,r}(s,t,z)=\sum^{r-1}_{l=0}P^l_{q,r}(z)\z_{q}(s-l,t-l,z).
\end{equation}
 For example, we have  
\begin{align*}
 \z_{q,2}(s,t,z)&=\z_{q}(s-1,t-1,z)-q[z-1]_q\z_{q}(s,t,z),\\
 \z_{q,3}(s,t,z)&=\frac{1}{1+q}\Bigl\{\z_{q}(s-2,t-2,z)\\
&\qquad
 -(q[z-1]_q+q^2[z-2])\z_{q}(s-1,t-1,z)+q^3[z-1]_q[z-2]_q\z_{q}(s,t,z)\Bigr\}. 
\end{align*}
 We now recall the analytic continuation of $\z_q(s,t,z)$ proved in 
 \cite{KawagoeWakayamaYamasaki}. Let $N\in\bN$. Put
 $f_q(x):=q^{xt}(1-q^{x+z})^{-s}$. Define the polynomial $b^{\e}_{j}(s)$
 $(0\le \e \le j)$ in $s$ by the following equation:
\[
 \frac{d^j}{dx^j}\{(1-q^{x+z})^{-s}\}
=(\log{q})^j\sum^{j}_{\e=0}b^{\e}_{j}(s)(1-q^{x+z})^{-s-\e}.
\]
 By the Leibniz rule, we have 
\[
 f\tp{k}_q(x)
=(\log{q})^kq^{xt}\sum^{k}_{\e=0}c^{\e}_{k}(s,t)(1-q^{x+z})^{-s-\e},
 \qquad c^{\e}_{k}(s,t):=\sum^{k}_{j=\e}\binom{k}{j}t^{k-j}b^{\e}_{j}(s).
\]
 Choosing $f(x)=f_q(x)$ and $M=N$ in \eqref{for:EM}, we have
\begin{multline}
\label{for:zq}
 \z_{q}(s,t,z)
=\frac{1}{2}\Bigl(\frac{1-q^z}{1-q}\Bigr)^{-s}-\sum^{N}_{k=1}\sum^{k}_{\e=0}\frac{B_{k+1}}{(k+1)!}c^{\e}_{k}(s,t)\Bigl(\frac{1-q^z}{1-q}\Bigr)^{-s-\e}\frac{(\log{q})^k}{(1-q)^{\e}}\\
+(1-q)^sI^0_{q,0}(s,t,z)+\frac{(-1)^N(\log{q})^{N+1}(1-q)^s}{(N+1)!}\sum^{N+1}_{\e=0}c^{\e}_{N+1}(s,t)I^{N+1}_{q,\e}(s,t,z),
\end{multline}
 where  
\[
 I^m_{q,\e}(s,t,z):=\int^{\infty}_{0}\tB_{m}(x)q^{xt}(1-q^{x+z})^{-s-\e}dx.
\]
 Note that $\tB_{0}(x)=1$. Recall now the Fourier expansion of
 $\tB_{m}(x)$ (see, e.g. \cite[p.\,191]{WhittakerWatson1927});  
\begin{equation}
\label{for:Fourier}
 \tB_{m}(x)=-m!\sumn\frac{e^{2\pi\I nx}}{(2\pi \I n)^m} \qquad
 (m\ge 2).
\end{equation}
 Put $u=q^{x+z}$. Then we have  
\begin{align}
 I^0_{q,0}(s,t,z)&=-\frac{q^{-zt}}{\log{q}}b_{q^z}(t,-s+1),
\label{for:I0}\\
 I^m_{q,\e}(s,t,z)&=\sumn\frac{m!e^{-2\pi\I nz}}{(2\pi\I
 n)^m}\frac{q^{-zt}}{\log{q}}b_{q^z}(\d n+t,-s-\e+1) \qquad (m\ge 2),
\label{for:Im} 
\end{align}
 where $\d=2\pi\I/\log{q}$. Here $b_w(\a,\b)$ is the incomplete beta
 function defined by the integral  
\[
 b_w(\a,\b):=\int^w_0u^{\a-1}(1-u)^{\b-1}du \qquad (0<\Re(w)<1).
\]
 This integral converges absolutely for $\Re(\a)>0$. Hence the function
 $b_w(\a,\b)$ is holomorphic for $\Re(\a)>0$ and for all
 $\b\in\bC$. Note that if $\Re(\b)>0$, we have
 $\lim_{w\to 1} b_w(\a,\b)=\B(\a,\b)$ where $\B(\a,\b)$ is the 
 beta function. Further, for any integer $N'\ge 2$, repeated use of
 integration by parts yields    
\begin{multline}
\label{for:ac-beta}
 \qquad  b_w(\a,\b)
=\sum^{N'-1}_{l=1}(-1)^{l-1}\frac{(1-\b)_{l-1}}{(\a)_l}w^{\a+l-1}(1-w)^{\b-l}\\   
+(-1)^{N'-1}\frac{(1-\b)_{N'-1}}{(\a)_{N'-1}}b_w(\a+N'-1,\b-N'+1).\qquad\qquad
\end{multline}
 As a function of $\a$, this expression gives an analytic continuation of
 $b_w(\a,\b)$ to the region $\Re(\a)>1-N'$. Hence the functions
 $I^0_{q,0}(s,t,z)$ and $I^m_{q,\e}(s,t,z)$ are meromorphically
 continued to the region $\Re(t)>1-N'$ for any integer $N'\ge 2$. Let
 $M\ge 0$ be an arbitrary large integer. Using the expressions
 \eqref{for:I0} and \eqref{for:Im}, and applying the formula
 \eqref{for:ac-beta} to $I^m_{q,\e}(s,t,z)$ with $N':=M-N+1\ge 2$, we
 see that the formula \eqref{for:zq} can be written as  
\begin{multline}
\label{for:EM-zq}
 \z_{q}(s,t,z)
=-\frac{q^{-zt}(1-q)^s}{\log{q}}b_{q^z}(t,-s+1)+\frac{1}{2}\Bigl(\frac{1-q^z}{1-q}\Bigr)^{-s}\\
+D^1_q(s,t,z;N,M)+D^2_q(s,t,z;N,M)+D^3_q(s,t,z;N,M),\qquad\qquad
\end{multline}
 where
\begin{align*}
 D^1_q(s,t,z;N,M):
&=-\sum^{N}_{k=1}\sum^{k}_{\e=0}\frac{B_{k+1}}{(k+1)!}c^{\e}_{k}(s,t)\Bigl(\frac{1-q^z}{1-q}\Bigr)^{-s-\e}\frac{(\log{q})^k}{(1-q)^{\e}},\\
 D^2_q(s,t,z;N,M):
&=\sum^{N+1}_{\e=0}\sum^{M-N}_{l=1}\sumn\frac{(-1)^{N+l-1}}{(2\pi\I
 n)^{N+1}}\frac{c^{\e}_{N+1}(s,t)(s+\e)_{l-1}q^{z(l-1)}}{(1-q)^l(\d
 n+t)_l}\\
&\qquad\qquad\qquad \times
 \Bigl(\frac{1-q^z}{1-q}\Bigr)^{-s-\e+1-l}\frac{(\log{q})^N}{(1-q)^{\e-1}},\\ 
 D^3_q(s,t,z;N,M)
:&=\sum^{N+1}_{\e=0}\sumn\frac{(-1)^{M+1}}{(2\pi\I
 n)^{N+1}}\frac{c^{\e}_{N+1}(s,t)(s+\e)_{M-N}q^{z(M-N)}}{(1-q)^{M-N}(\d n+t)_{M-N}}\frac{(\log{q})^{N+1}}{(1-q)^{\e}}\\
&\qquad\qquad\qquad \times \int^{\infty}_{0}e^{2\pi\I nx}q^{x(t+M-N)}\Bigl(\frac{1-q^{x+z}}{1-q}\Bigr)^{-s-\e-M+N}dx.
\end{align*}
 The equation \eqref{for:EM-zq} gives an analytic continuation of
 $\z_q(s,t,z)$ to the region $\Re(t)>1-N'=N-M$. Note that, by the fact
 $c^{k}_{k}(s,t)=(s)_k$ and \eqref{for:Fourier} again, we have  
\begin{align}
 \limup D^1_q(s,t,z;N,M)
&=\sum^N_{k=1}\frac{B_{k+1}}{(k+1)!}(s)_kz^{-s-k},
\label{for:limD1}\\
 \limup D^2_q(s,t,z;N,M)
&=\sum^{M}_{l=N+1}\frac{B_{l+1}}{(l+1)!}(s)_lz^{-s-l},
\label{for:limD2}\\
 \limup D^3_q(s,t,z;N,M)
&=-\frac{(s)_{M+1}}{(M+1)!}\int^{\infty}_{0}\widetilde{B}_{M+1}(x)(x+z)^{-s-M-1}dx.\label{for:limD3}  
\end{align}
 Therefore, by \eqref{for:zqr-zq} and \eqref{for:EM-zq}, we obtain the
 following  

\begin{prop}
\label{prop:zqr}
 For any integers $N_l\ge 1$ and $M_l\ge N_l+1$ $(0\le l \le r-1)$, we
 have   
\begin{align*}
 \z_{q,r}(s,t,z)
&=-\frac{(1-q)^{s-(r-1)}}{\log{q}}\sum^{r-1}_{l=0}P^l_{q,r}(z)q^{-z(t-l)}(1-q)^{r-1-l}b_{q^z}(t-l,-s+l+1)\\
&\qquad +\frac{1}{2}\sum^{r-1}_{l=0}P^l_{q,r}(z)\Bigl(\frac{1-q^z}{1-q}\Bigr)^{-s+l}+\sum^{r-1}_{l=0}P^l_{q,r}(z)D^1_q(s-l,t-l,z;N_l,M_l)\\
&\qquad +\sum^{r-1}_{l=0}P^l_{q,r}(z)D^2_q(s-l,t-l,z;N_l,M_l)+\sum^{r-1}_{l=0}P^l_{q,r}(z)D^3_q(s-l,t-l,z;N_l,M_l).
\end{align*}
 This gives an analytic continuation of $\z_{q,r}(s,t,z)$ to the region
 $\Re(t)>M'$ where $M':=\max\{N_l-M_l+l\,|\,0\le l\le r-1\}$.\qed 
\end{prop}

\subsection{Proof of Theorem\,\ref{thm:main}}

 Note the following lemma.  

\begin{lem}
\label{lem:P}
 It holds that 
\[
 \limup P^l_{q,r}(z)=P^l_r(z).
\]
\end{lem}
\begin{proof}
 By the definition of $P^l_{q,r}(z)$, it is sufficient to show
\begin{equation}
\label{for:suff}
 (-1)^{r-1-l}\sum_{1\le m_1<\cdots<m_{r-1-l\le r-1}}(z-m_1)\cdots
 (z-m_{r-1-l})
=\sum^{r-1}_{k=l}\binom{j}{l}s(r,j+1)(-z)^{j-l}.
\end{equation}
 Notice that the left hand side of \eqref{for:suff} is equal to the
 coefficient of $x^l$ in the polynomial 
 $p_r(x):=\prod^{r-1}_{j=1}(x-(z-j))$ in $x$. Since
 $p_r(x)=(x-z)_r/(x-z)$, we have    
\[
 p_r(x)=\sum^{r-1}_{j=0}s(r,j+1)(x-z)^{j}
=\sum^{r-1}_{l=0}\Bigl(\sum^{r-1}_{j=l}\binom{j}{l}s(r,j+1)(-z)^{j-l}\Bigr)x^l.
\]
 Hence the desired formula follows.\qed
\end{proof}

 We are ready to prove the main theorem.

\begin{Proof}{\Thm{main}}
 We first show the sufficiency. Let $t=s-\n$ for $\n\in\bN$. Notice
 that, by \eqref{for:zqr-zq}, we have
 $\z\nn_{q,r}(s,z)=\sum^{r-1}_{l=0}P^l_{q,r}(z)\z\nn_q(s-l,z)$. Hence,
 by \cite[Theorem\,$2.1.$]{KawagoeWakayamaYamasaki}, Lemma\,\ref{lem:P}
 and \eqref{for:zr-z}, we have
\[
 \limup \z\nn_{q,r}(s,z)
=\sum^{r-1}_{l=0}P^l_{r}(z)\z(s-l,z)=\z_r(s,z) \qquad (s\in\bC).
\] 

 We next show the necessity. Suppose that $\limup \z_{q,r}(s,t,z)$
 exists and satisfies $\limup \z_{q,r}(s,t,z)=\z_r(s,z)$ for all
 $s\in\bC$ with some meromorphic function $t=\vp(s)$. Then, by
 \Prop{zr}, \Prop{zqr}, \Lem{P}, \eqref{for:limD1}, \eqref{for:limD2}
 and \eqref{for:limD3}, it is necessary to hold
\[
 -\limup \frac{(1-q)^{s-(r-1)}}{\log{q}}\sum^{r-1}_{l=0}P^l_{q,r}(z)q^{-z(t-l)}(1-q)^{r-1-l}b_{q^z}(t-l,-s+l+1)=\sum^{r-1}_{l=0}\frac{P^l_r(z)}{s-l-1}z^{-s+l+1}.
\]
 Assume $\Re(s)<1$. Since $\limup (1-q)^{s-(r-1)}/\log{q}$ diverges, it
 is necessary to hold    
\begin{equation}
\label{for:sufff}
 \limup \sum^{r-1}_{l=0}P^l_{q,r}(z)q^{-z(t-l)}(1-q)^{r-1-l}b_{q^z}(t-l,-s+l+1)=0.
\end{equation}
 Notice that $\limup b_{q^z}(t-l,-s+l+1)=\B(t-l,-s+l+1)$
 for all $l$ $(0\le l\le r-1)$. 
 Further, since the left hand side of \eqref{for:sufff} is equal to 
 $\B(t-r+1,-s+r)=\G(t-r+1)\G(-s+r)/\G(t-s+1)$, we have
 $t-s+1\in\bZ_{\le 0}$, whence $t=\vp(s)=s-\n$ for some positive integer
 $\n\in\bN$ in the region $\Re(s)<1$. Since $\vp(s)$ is meromorphic on
 $\bC$, we have $\vp(s)=s-\n$ for all $s\in\bC$. This proves the
 theorem.\qed  
\end{Proof}

\section{Remarks on $q$-analogues of $\z_r(s,z;\bsym{\om})$}
\label{sec:general}

 We introduce here a $q$-analogue of the Barnes multiple zeta function
 $\z_r(s,z;\bsym{\om})$ for a general parameter
 $\bsym{\om}:=(\om_1,\ldots,\om_r)$. Assume $\om_i>0$ $(1\le i \le r)$
 and $\Re(z)>0$. We define a $q$-analogue of
 $\z_{q,r}(s,t,z;\bsym{\om})$ by the series     
\[
 \z_{q,r}(s,t,z;\bsym{\om}):
=\sum_{n_1,\ldots,n_r\ge 0}\frac{q^{n_1\om_1t+n_2\om_2(t-1)+\cdots+n_r\om_r(t-r+1)}}{{[n_1\om_1+\cdots+n_r\om_r+z]_q}^s}.
\]
 We put $\z\nn_{q,r}(s,z;\bsym{\om}):=\z_{q,r}(s,s-\n,z;\bsym{\om})$ for 
 $\n\in\bN$. The series $\z_{q,r}(s,t,z;\bsym{\om})$ converges
 absolutely for $\Re(t)>r-1$. It is clear that
 $\z_{q,r}(s,t,z)=\z_{q,r}(s,t,z;{\bf 1}_r)$ where
 ${\bf 1}_r:=(\underbrace{1,1,\ldots,1}_r)$. By the following
 proposition, $\z_{q,r}(s,t,z;\bsym{\om})$ is continued meromorphically
 to the whole $s$, $t$-plane. The proof can be obtained by the similar
 way to \cite[Proposition\,1]{KanekoKurokawaWakayama2003} and
 \cite[Proposition\,2.9]{KawagoeWakayamaYamasaki}.    

\begin{prop}
\label{prop:ac-Gzqr}
 $(i)$ The function $\z_{q,r}(s,t,z;\bsym{\om})$ can be written as 
\begin{equation}
\label{for:ac-Gzqr}
 \z_{q,r}(s,t,z;\bsym{\om})
=(1-q)^s\sum^{\infty}_{l=0}\binom{s+l-1}{l}q^{lz}\prod^r_{j=1}(1-q^{\om_j(t-j+1+l)})^{-1}.
\end{equation}
 This gives a meromorphic continuation of $\z_{q,r}(s,t,z;\bsym{\om})$ to the
 whole $s$, $t$-plane with simple poles at
 $t\in j-1+\bZ_{\le 0}+\d_j\bZ$\ $(1\le j \le r)$. Here
 $\d_j:=2\pi\I/(\om_j\log{q})$. 

 $(ii)$ The function $\z\nn_{q,r}(s,z;\bsym{\om})$ can be written as 
\begin{equation}
\label{for:ac-Gzqrn}
 \z\nn_{q,r}(s,z;\bsym{\om})
=(1-q)^s\sum^{\infty}_{l=0}\binom{s+l-1}{l}q^{lz}\prod^r_{j=1}(1-q^{\om_j(s-\n-j+1+l)})^{-1}.
\end{equation}
 This gives a meromorphic continuation of $\z\nn_{q,r}(s,z;\bsym{\om})$
 to the whole plane $\bC$ with simple poles at the points in 
\[
 \begin{cases}
  \ j+\d_i\bZ\bslo & (\,j\in\bZ_{\le 0},\ 1\le i \le r),\\
  \ j+\d_i\bZ & (\,1\le j \le \n,\ 1\le i \le r),\\
  \ \n+j+\d_i\bZ & (\,1\le j \le r-1,\ j+1 \le i \le r).
 \end{cases}
\]
 In particular, the poles of $\z\nn_{q,r}(s,z;\bsym{\om})$ on the real
 axis are given by $s=1,2,\ldots,r,r+1,\ldots,r+\n-1$.

$(iii)$\ Let $m\in\bZ_{\ge 0}$. Then we have 
\begin{multline}
\label{for:specialvelue}
 \z\nn_{q,r}(-m,z;\bsym{\om})
=(1-q)^{-m}\Biggl\{\sum^{m}_{l=0}(-1)^l\binom{m}{l}q^{lz}\prod^r_{j=1}(1-q^{\om_j(-m-\n+l-j+1)})^{-1}\\
+\frac{q^{(m+\nu-1)z}}{\log{q}}\sum^r_{l=1}\frac{(-1)^{m+1}m!(l+\n-2)!\,q^{lz}}{(l+m+\n-1)!\,\om_l}\prod^r_{j=1 \atop j\ne l}(1-q^{\om_j(l-j)})^{-1}\Biggr\}.
\end{multline}
\end{prop}
\begin{proof}
 The formula \eqref{for:ac-Gzqr} is obtained by the binomial theorem,
 whence \eqref{for:ac-Gzqrn} immediately follows. The formula
 \eqref{for:specialvelue} is derived from the fact
 $(s+m)/(1-q^{\om_l(s+m)})=-1/(\om_l\log{q})+O(s+m)$ as $s\to -m$.\qed  
\end{proof}

\noindent
 Theses facts motivate the

\begin{conj}
\label{conj:general}
 Let $t=\vp(s)$ be a meromorphic function on $\bC$. Then the formula 
\[
 \limup \z_{q,r}(s,\vp(s),z;\bsym{\om})=\z_r(s,z;\bsym{\om}) \qquad (s\in\bC)
\]
 holds if and only if the function $\vp(s)$ can be written as
 $\vp(s)=s-\n$ for some $\n\in\bN$. 
\end{conj}
 
\noindent
 In fact, since $\z_1(s,z;\bsym{\om})={\om}^{-s}\z(s,z/{\om})$ and
 $\z_{q,1}(s,t,z;\bsym{\om})={[\om]_q}^{-s}\z_{q^{\om}}(s,t,z/{\om})$ for
 $\om>0$, \Conj{general} is true for $r=1$ by \eqref{for:EM-z} and
 \eqref{for:EM-zq}.

\appendix

\section{Associated $q$-analogue of the gamma function}
\label{appen:gamma}

 In this appendix, we introduce a $q$-analogue of the gamma
 function defined via the $q$-analogue of the Hurwitz zeta function:
\[
 \tz_q(s,z):
=\z^{(1)}_q(s,z)\times q^{z(s-1)}
=\sum^{\infty}_{n=0}\frac{q^{(n+z)(s-1)}}{{[n+z]_q}^s} \qquad 
 (\Re(s)>1).
\]
 Note that by \eqref{for:ladder zqr}, we have 
\begin{equation}
\label{for:ladder tzq}
 \tz_q(s,z)=\tz_q(s,z+1)+\frac{q^{z(s-1)}}{{[z]_q}^s}.
\end{equation}
 Imitating the Lerch formula \cite{Lerch1894} (the zeta regularization)
\[
 \frac{\p}{\p s}\z(s,z)\Bigl|_{s=0}=\log\frac{\G(z)}{\sqrt{2\pi}},
\]
 we define a $q$-analogue $\tG_q(z)$ of the gamma function by    
\[
 \tG_q(z):
=\exp\Bigl(\frac{\p}{\p s}\tz_q(s,z)\Bigl|_{s=0}
-\frac{\p}{\p s}\tz_q(s,1)\Bigr|_{s=0}\Bigr).
\] 
 Then the function $\tG_q(z)$ is well-defined as a single valued
 meromorphic function. Indeed, let 
\[
 \tz_q(s,z)=a_0(z;q)+a_1(z;q)s+a_2(z;q)s^2+\cdots
\]
 be the Taylor expansion of $\tz_q(s,z)$ around $s=0$. Note that
 $\tz_q(s,z)$ is holomorphic at $s=0$. Assume $\Re(z)>0$. Then, by
 \Prop{ac-Gzqr}, $\tz_q(s,z)$ has the following expression; 
\begin{equation}
\label{for:ac-tzq}
 \tz_q(s,z)
=(1-q)^s\sum^{\infty}_{n=0}\binom{s+n-1}{n}\frac{q^{z(s-1+n)}}{1-q^{s-1+n}}.
\end{equation}
 Hence one can calculate the coefficient $a_1(z;q)$ by the same manner
 performed in \cite{KurokawaWakayama2003} as 
\begin{equation}
\label{def:diff0}
 a_1(z;q)
=\sum^{\infty}_{n=2}\frac{1}{n}\frac{q^{(n-1)z}}{1-q^{n-1}}-z+\frac{1}{2}
+\frac{1-z(1-q)}{(1-q)^2}q^{1-z}\log{q}-\Bigl(\frac{q^{1-z}}{1-q}+\frac{1}{\log{q}}\Bigr)\log{(1-q)}.
\end{equation}
 Therefore $\tG_q(z)$ is meromorphic in the region $\Re(z)>0$. If
 $-1<\Re(z)<0$, by the ladder relation \eqref{for:ladder tzq}, we have  
\[
 a_1(z;q)=q^{-z}\log{q^z}-q^{-z}\log\Bigl(\frac{1-q^z}{1-q}\Bigr)+a_1(z+1;q).
\] 
 Hence we have
\begin{equation}
\label{for:nori}
 \tG_q(z)=(q^{-z}[z]_q)^{-q^{-z}}\tG_q(z+1).
\end{equation}
 This gives a meromorphic continuation of $\tG_q(z)$ to the region
 $\Re(z)>-1$. Repeating the same procedure, we see that $\tG_q(z)$ can
 be extended as a meromorphic function on $\bC$. 

 From \Thm{main}, by the Lerch formula, we have immediately 
\begin{align}
\label{for:lim-q-gamma}
 \limup \tG_q(z)=\G(z) \qquad (z\notin -\bN_0). 
\end{align}
 Moreover, $\tG_q(z)$ satisfies the following properties.

\begin{prop}
\label{prop:q-Bohr}
 We have
\begin{align}
 &\tG_q(z+1)=(q^{-z}[z]_q)^{q^{-z}}\tG_q(z),
\label{for:ladder tG_q}\\
&\tG_q(1)=1,
\label{for:tG_q1}\\
 &\frac{d^2}{dz^2}\log{\tG_q(z+1)}\ge 0 \qquad (z\ge 0).
\label{for:log conv}\qquad\qquad
\end{align} 
 In particular, for a positive integer $n$, we have
\begin{equation}
\label{for:tGqn}
 \tG_q(n+1)=q^{-\sum^{n}_{k=1}kq^{-k}}\prod^{n}_{k=1}([k]_q)^{q^{-k}}.
\end{equation}
\end{prop}
\begin{proof}
 By the definition of $\tG_q(z)$, \eqref{for:tG_q1} is obvious. The
 formula \eqref{for:ladder tG_q} is clear from \eqref{for:nori}. The
 assertion \eqref{for:tGqn} follows from \eqref{for:ladder tG_q} and
 \eqref{for:tG_q1} by induction. To show the inequality \eqref{for:log
 conv}, take the 
 logarithm of $\tG_q(z)$:   
\begin{multline}
\label{for:log-tGq}
 \quad \log\tG_q(z)
=\sum^{\infty}_{n=2}\frac{1}{n}\frac{q^{z(n-1)}-q^{n-1}}{1-q^{n-1}}-z+1\\
+\frac{q^{-z}\bigl(1-(1-q)z\bigr)-1}{(1-q)^2}q\log{q}+\frac{1-q^{1-z}}{1-q}\log{(1-q)}.\qquad\qquad
\end{multline}
 We calculate as  
\[
 \frac{d^2}{dz^2}\log{\tG_q(z+1)}
=(\log{q})^2\sum^{\infty}_{n=2}\frac{(n-1)^2}{n}\frac{q^{(z+1)(n-1)}}{1-q^{n-1}}
+\frac{(\log{q})^2q^{-z}}{(1-q)^2}\eta_q(z),
\]
 where
 $\eta_q(z):=(\log{q})\bigl(1-(1-q)(z+1)\bigr)-(1-q)\log{(1-q)}+2(1-q)$.
 Therefore, it suffices to show that $\eta_q(z)\ge 0$ for all $0<q<1$ if 
 $z\ge 0$, and this is indeed true. In fact, since
 $\frac{d}{dq}\eta_q(z)\le 0$ for $0<q<1$, we conclude that
 $\eta_q(z)\ge \limup \eta_q(z)=0$. Hence the proposition follows.\qed   
\end{proof}

\begin{remark}
 One can find the similar formulas to \eqref{for:ladder tG_q},
 \eqref{for:tG_q1} and \eqref{for:log conv} in the $q$-analogue of
 the Bohr-Morellup theorem for the Jackson $q$-gamma function in
 \cite{Askey1978}. It has not yet been clarified that theses properties
 characterize the function $\tG_q(z)$. 
\end{remark}

\smallbreak

 By the expression \eqref{for:ac-tzq} again, $\tz_q(s,z)$ has the
 following Laurent expansion around $s=1$:   
\begin{equation}
\label{for:s1}
 \tz_q(s,z)=\frac{q-1}{\log{q}}\frac{1}{s-1}+\g_q(z)+O(s-1) \qquad (\Re(z)>0),
\end{equation} 
 where
\begin{equation}
\label{def:gq}
 \g_q(z)
:=\sum^{\infty}_{n=1}\frac{q^{nz}}{[n]_q}+(1-q)\Bigl(-z+\frac{1}{2}-\frac{\log{(1-q)}}{\log{q}}\Bigr).
\end{equation}
 We next show a $q$-analogue of the Lerch limit formula \cite{Lerch1894}:  
\begin{equation}
\label{for:Lerch}
 \lim_{s\to 1}\Bigl(\z(s,z)-\frac{1}{s-1}\Bigr)=-\frac{\G'}{\G}(z).
\end{equation}

\begin{prop}
\label{prop:q-Lerch}
 It holds that 
\begin{equation}
\label{for:q-Lerch}
 \g_q(z)
=\lim_{s\to 1}\Bigl(\tz_q(s,z)-\frac{q-1}{\log{q}}\frac{1}{s-1}\Bigr)
=-\frac{q-1}{\log{q}}\frac{\tG'_q}{\tG_q}(z)+C_q(z) \qquad (\Re(z)>0),
\end{equation}
 where 
\begin{multline*}
\qquad C_q(z):
=\sum^{\infty}_{n=1}\frac{1}{n+1}\frac{q^{nz}}{[n]_q}-q^{1-z}-\frac{\log{q}}{1-q}(1-(1-q)z)q^{1-z}+\frac{1-q}{\log{q}}\\
-q^{1-z}\log{(1-q)}-\frac{1-q}{\log{q}}\log{(1-q)}+\Bigl(-z+\frac{1}{2}\Bigr)(1-q)\qquad\qquad
\end{multline*}
 and $\limup C_q(z)=0$. Put $\g_q:=\g_q(1)$. Then we have, in
 particular, $\limup \g_q=\g$ where $\g=0.577215\ldots$ denotes the Euler 
 constant. 
\end{prop}
\begin{proof}
 By \eqref{for:log-tGq}, we have 
\begin{multline}
\label{for:log-deriv-Gq}
 \qquad\frac{\tG_q'}{\tG_q}(z)
=\frac{\log{q}}{1-q}\Bigl(\sum^{\infty}_{n=1}\frac{q^{zn}}{[n]_q}-\sum^{\infty}_{n=1}\frac{1}{n+1}\frac{q^{nz}}{[n]_q}\Bigr)-1\\
+\frac{(1-q)+\bigl(1-(1-q)z\bigr)\log{q}}{(1-q)^2}q^{1-z}\log{q}+\frac{\log{q}}{1-q}q^{1-z}\log{(1-q)}.\qquad
\end{multline}
 Plugging \eqref{def:gq} into \eqref{for:log-deriv-Gq}, we obtain the
 formula \eqref{for:q-Lerch}. It is straightforward to show the fact
 $\limup C_q(z)=0$ when $\Re(z)>0$. Hence we have $\limup \g_q=\g$ by
 the limit formulas \eqref{for:lim-q-gamma}, \eqref{for:Lerch} and the
 facts $\G(1)=1$, $\G'(1)=-\g$. This completes the proof.\qed   
\end{proof}

\begin{remark}
 The $q$-analogue of the Lerch limit formula obtained in this paper is
 different from the one given in \cite{KurokawaWakayama2003}.
\end{remark}
\smallbreak

 As a final remark, we give a $q$-analogue of the Gauss-Legendre formula.

\begin{prop}
\label{prop:Gauss-Legendre}
 Let $N\in\bN$. Then we have 
\[
 [N]^{[1-Nz]_q}_q\tG_{q^N}\bigl(\frac{1}{N}\bigr)\cdots\tG_{q^N}\bigl(\frac{N-1}{N}\bigr)\tG_q(Nz)
=\tG_{q^N}(z)\tG_{q^N}\bigl(z+\frac{1}{N}\bigr)\cdots\tG_{q^N}\bigl(z+\frac{N-1}{N}\bigr).
\]
\end{prop}
\begin{proof}
 The proof is straightforward from \eqref{for:log-tGq}.\qed 
\end{proof}


\end{document}